\documentstyle[12pt]{article}

\advance\hoffset by -7mm
\makeatletter
\@addtoreset{equation}{section}
\makeatother

\input amssym.def  %% use .def & .tex for \Bbb commands%%
\input amssym.tex

\renewcommand{\title}[1]{\null\vspace{25mm}

\noindent{\Large{\bf #1}}\vspace{10mm}

\noindent {\large By }}
\newcommand{\authors}[1]{\noindent{\large #1}\vspace{3mm}

}
\newcommand{\address}[1]{\noindent #1\vspace{5mm}

}
\renewcommand{\abstract}[1]{\vspace{19mm}

\noindent{\small{\em Abstract.} #1}\vspace{2mm}

}

\begin{document}

\title{Variational properties of a nonlinear elliptic equation and
  rigidity}
\authors{M L Bialy}
\address{School of Math Sciences, Tel-Aviv University, Tel-Aviv 69978,
  Israel\\
Email:  bialy@math.tau.ac.il}
\authors{and R S MacKay}
\address{Nonlinear Centre, DAMTP, Silver Street, Cambridge CB3 9EW,
  UK\\
Email:  R.S.MacKay@damtp.cam.ac.uk}

\abstract{
We consider in this paper elliptic
equations which are perturbations of Laplace's equation by a compactly
supported potential.  We show that in dimension greater than three for a
wide class of potentials all the solutions are globally minimising.
However, in dimension two the situation is different.  We show that
for radially symmetric potentials there always exist solutions which
are not locally minimal unless the potential vanishes identically.  We
discuss the relations of this with the so-called Hopf rigidity
phenomenon.
}

\section{Introduction and main results}

In this paper we discuss variational properties of classical solutions
of the nonlinear equation
\begin{equation}
\Delta u = -V'_u \left( u, x_1, \dots, x_n \right)
\label{11}
\end{equation}
which is the Euler-Lagrange equation of the functional
\begin{equation}
I(u)= \int {1 \over 2} \left( \nabla u \right)^2 - V \left( u, x_1,
  \dots, x_n \right) dx_1 \dots dx_n
\label{12}
\end{equation}

The main question which is addressed here is the following:  under
what conditions on the potential $V$ are all classical solutions of
(\ref{11}) globally minimising for the functional (\ref{12})?
By a global minimiser we mean a smooth function on a domain of $\Bbb R^n$
minimising the integral (1.2) over all bounded subdomains with
smooth boundaries with respect to smooth functions with the same
boundary values.

The motivation for this question comes from variational problems of
classical mechanics.  In geodesic problems it is well known that all
geodesics are globally minimising on manifolds of negative sectional
curvature.  However, it was shown first by E Hopf and L Green that the
situation is completely different for Riemannian (two) tori or for
Riemannian planes which are flat outside a compact set.  They proved
(\cite{[Ho]}, \cite{[G]}) that for these manifolds there always exists
geodesics with conjugate points and therefore non-minimal, unless the
metric is flat.  We refer the reader to \cite{[B-I]}, \cite{[C]} for
higher dimensional generalisations of E Hopf and Green's theorems and to
\cite{[C-K]}, \cite{[C-F]}, \cite{[K]}, \cite{[He]} for very important
previous developments.

It was observed first in \cite{[B]} that the Hopf phenomenon is not
entirely
Riemannian.  In \cite{[B-P]} we have shown that a similar type of
rigidity holds true for Newton's equations with periodic or compactly
supported potentials.

For the equation (\ref{11}) with periodic potentials it was shown in
\cite{[M1]}, \cite{[M2]} that far-going generalisations of Aubry-Mather
and KAM theories apply.  Using these theories one can construct
families of minimal solutions which form laminations or sometimes even
foliations of the configuration torus.
It is an interesting open question, however, if it
happens that for all ``slope'' vectors these laminations are genuine
foliations.  This question is tightly related to the one we are addressing
in this paper, since all the leaves of a foliation are globally minimal.

We shall assume throughout this paper that the potential $V$ is
compactly supported. The main reason for this is to make the analogy
with the above mentioned situations of Hopf rigidity
and to exclude, in particular, the case with $V''_{uu} \leq 0$
everywhere which is analogous to non-positive curvature.
Nevertheless, several of our results apply in more generality,
as we shall indicate.

We will show that in case of dimensions greater than 2 there are
many potentials such that all solutions of (\ref{11}) are minimising
(Theorem 1).
In dimension 2, however, at least for radially symmetric potentials
there always exist non-minimal solutions of (\ref{11}) unless $V$
vanishes identically (Theorem 2). We state theorems 1 and 2 here.

%Theorem 1.3
\newtheorem{theorem}{Theorem}
\begin{theorem}
\label{T13}
For $n \geq 3$ let $V$ be a compactly supported potential on
$\Bbb R^{n+1}$.  Assume that $V''_{uu} (u,x) \leq U(x)$
for some function $U$ such that either

$(A) \quad U(x) \leq \left( {n-2 \over 2} \right)^2 {1 \over ||x-x_0||^2}$
for some point $x_0 \in \Bbb R^n$, or

$(B) \quad ||U||_{n/2} \leq {n(n-2) \over 4} |\Bbb S^n|$
where $|\Bbb S^n|$ is the volume of the unit n-sphere.

\noindent Then {\it any} solution of
(\ref{11}) is globally minimising.
\end{theorem}

%Theorem 1.5
\begin{theorem}
\label{T15}
Let $V(u, x_1, x_2)$ be a radially symmetric compactly supported
potential $(n=2)$.  There always exist radial non-minimal solutions of
(\ref{11}) unless $V$ vanishes identically.
\end{theorem}

Our approach to the proof of the theorem 2 is based on
the reduction to a Newton equation with the potential supported
in the semi-strip. It turns out that a technique analogous to Hopf's
can be applied for such a shape of support.

In Section 3 we will consider the case of radially symmetric
potentials for $n\geq 3$.  We illustrate the minimality property of
radial solutions of (\ref{11}) organising them in foliations (Theorem 3).

\noindent{\large \bf Acknowledgements}

\noindent
This paper was written while the first author visited the University
of Cambridge.  We would like to thank the EPSRC for their support.

%Theorem 1.3 (section 2)
\section{Proof of Theorem {\ref{T13}}}

Let $u$
be a solution of (\ref{11}) in a domain $\Omega$ and $\upsilon$ be any
function with the same boundary values.  Set $\xi = \upsilon -u$ and
write
\begin{equation}
I(\upsilon) - I(u) = \int_\Omega \nabla u \nabla \xi + {(\nabla \xi)^2
  \over 2} - \left( V(u + \xi, x) - V(u,x) \right) d^nx
\label{31}
\end{equation}
Integrating the first term by parts and using the equation (\ref{11}) we have
$$
\int_\Omega \nabla u \nabla \xi dx = \int_\Omega V'_u (u,x) \xi \, d^nx
$$
By the assumption on $V$ we have
$$
V(u+\xi, x) - V(u,x) - V'_u(u,x)\xi \leq {1 \over 2} U(x) \xi^2
$$
Substituting this in (\ref{31}) we obtain
\begin{equation}
I(\upsilon) - I(u) \geq \int_\Omega {(\nabla \xi)^2\over 2} - {1 \over
  2} U(x) \xi^2 \, d^nx
\label{32}
\end{equation}

We show that under hypothesis (A) or (B)
the last integral is always positive unless $\xi \equiv
0$.  Indeed, for case (A) introduce spherical coordinates centred at $x_0$,
$r=||x-x_0||$.  The required assertion follows from the following

%Lemma 3.3
\newtheorem{lemma}{Lemma}
\begin{lemma}
\label{L33}
For any function $\xi(r)$ defined on $\left[r_1, r_2\right] \left( 0
  \leq r_1 < r_2 \leq \infty \right)$ with $\xi \left( r_2 \right) =0$
  it follows that
$$
\int^{r_2}_{r_1} r^{n-1} \left( \left( \xi'\right)^2 - \left({n-2
      \over 2} \right)^2 {\xi^2 \over r^2} \right) dr \geq 0
$$
with equality for $\xi \equiv 0$ only.
\end{lemma}

\noindent{\bf Proof of Lemma {\ref{L33}}}

Introduce the new function $\varphi (r) = \xi (r) \cdot r^{n/2-1}$ and
substitute into the integral.  We have
\begin{eqnarray*}
&&\int^{r_2}_{r_1} r^{n-1} \left[ \left( r^{1 - {n \over 2}} \varphi' +
    \left({2-n \over 2} \right) r^{-{n \over 2}}\varphi \right)^2 -
  {(n-2)^2 \over 4} r^{-n} \varphi^2 \right] dr\\
&&= \int^{r_2}_{r_1} \left( r \cdot \left( \varphi' \right)^2 + (2-n)
    \varphi \varphi' \right) dr\\
&&= {n-2 \over 2} \varphi^2 \left( r_1 \right) + \int^{r_2}_{r_1} r
    \left( \varphi' \right)^2 dr
\end{eqnarray*}
Since $n\geq 3$ we see that the last expression is always non-negative
and equals zero only for $\xi \equiv 0$.  This completes the proof of
Lemma {\ref{L33}} and of the theorem in case (A).

Alternatively, under hypothesis (B), Sobolev's inequality (e.g.~\S8.3 of
\cite{[LL]}) implies that
$$\int |\nabla \xi|^2 \, d^{n}x \ge S_n ||\xi||_{2n\over n-2}^2,$$
where $$S_n = {n(n-2)\over 4}|\Bbb S^n|.$$
Then H\"older's inequality applied to the second term of the integrand
(as in \S11.3 of \cite{[LL]}) yields
$$|\int U \xi^2 \, d^{n}x| \le ||\xi^2||_{n\over n-2} ||U||_{n\over 2} .$$
Thus
$$ I(v) - I(u) \ge {1\over 2} (S_n - ||U||_{n\over 2}) ||\xi||_{2n\over
n-2}^2 .$$
This completes the proof in case (B).
$\Box$

\newtheorem{remark}{Remark}
\begin{remark}
Lemma 1 and its proof have a precursor in \cite{[CH]} (Ch6 \S 5),
and the whole of Theorem 1 fits in the general domain of ``absence of
bound states" described for example in \cite{[RS]} (\S 8.3).
\end{remark}
% Remark 3.4
\begin{remark}
\label{R34}
Our proof shows in addition that any solution is a
non\--degenerate minimum for (\ref{12}).
In a different way one can say this as follows:  for any solution $u$
on $\Omega$, the linearised equation
$$
\Delta \xi + V''_{uu} (u,x) \xi = 0
$$
has only the trivial solution satisfying the zero boundary conditions $\xi
|_{ \partial \Omega} = 0$.
\end{remark}

\section{Radial potentials}
Let us consider now the case of radial, compactly supported potentials
$V(u,r)$ satisfying $V''_{uu} (u,r) \leq {1 \over 4r^2}$.  Radial
solutions $u(r)$ of (\ref{11}) are given by the following:
$$
u'' + {n-1 \over r} u' + V'_u (u,r) = 0
$$
It follows from Theorem 1A that for $n \geq 3$ all radial
solutions are without
conjugate points.  Moreover, Lemma {\ref{L33}} gives that the points
$r_1$ and $\infty$ are not conjugate for any $r_1>0$.
This fact enables us to organise the solutions into foliations in
the following way.
Denote by $N_A$ the class of all those solutions which can be written
as $u(r, \alpha) = {\alpha \over r^{n-2}} + A$ for some $\alpha$
outside the support of $V$.

%Theorem 3.5
\begin{theorem}
\label{T35}
For compactly supported radial potential $V$ satisfying $V''_{uu} (u,r)
\leq {1 \over 4r^2}$ the set $N_A$ is totally ordered and the graphs of
the solutions define a smooth foliation of $\Bbb R^{n+1}- \Bbb R (u) \times
\{ x=0 \}$.
\end{theorem}

\noindent{\bf Proof of the theorem}

Given $A$ we have to show that the function $u(r,\alpha)$ is monotone
in $\alpha$.  Indeed the function $\xi (r) = {\partial u \over
  \partial \alpha}$ is a solution of the linearised equation.  Note
that by definition $\xi (\infty) =0$.  Then the non-conjugacy property
implies that $\xi>0$ and then it is easy to complete the proof. $\Box$

% Remark 3.6
\begin{remark}
\label{R36}
In some cases there is another way of organising the set of all
solutions into foliations.  Assume we are given a radial potential
with $V(u,r) \equiv 0$ for $0<r \leq r_0$ and $r \geq R_0$ satisfying
the inequality of Theorem {\ref{T35}}.  Define the set $M_A$ of all
those solutions which can be written as $u\left( r, \alpha \right) =
{A \over r^{n-2}}+ \alpha$ when $r \leq r_0$.  Then one shows that for
any given $A$ the set $M_A$ is ordered and the graphs of the solutions
smoothly foliate the space $\Bbb R^{n+1} - \Bbb R(u) \times \{ x=0\}$.  In
addition $M_0$ foliates the whole of $\Bbb R^{n+1}$.  In order to check the
order property one shows that the linearised equation has no focal
points in the following sense:  any solution satisfying $\dot{\xi}
\left( r_1 \right) = \xi \left( r_2 \right) = 0$ is trivial provided
$r_1 < r_2$ (this is not necessarily true for $r_1>r_2$).  Then one
proceeds exactly as in the proof of theorem {\ref{T35}}.
\end{remark}

\section{Rigidity for the case $n=2$}

Let $u(r)$ be a radial solution of the equation (\ref{11}) for
compactly supported rotationally symmetric $V(u,r)$, $V(u,r) \equiv 0$
for $|u|\geq U$ or $r>R$.  With the substitution $r=e^t$ the equation
for $u$ as a function of $t$ can be written in the
Hamiltonian form
\begin{eqnarray}
\dot u &=& p \nonumber \\
\dot p &=& -e^{2t} W'_u(u,t)
\label{41}
\end{eqnarray}
The Hamiltonian function of (\ref{41}) is
$$
H={1 \over 2} p^2 + e^{2t} W(u,t)
$$
with the function $W(u,t)=V \left( u, e^t \right)$.
Note that the support of $W$ is contained in the semi-strip $\Pi = \{
|u| \leq U, t \leq T = ln R \}$.

We shall prove the following rigidity result which implies theorem 2.

%Theorem 4.2
\begin{theorem}
There always exist solutions with conjugate points for (\ref{41})
unless $W$ vanishes identically.
\end{theorem}

The strategy of the proof will follow the original one of E Hopf, but
will take special care about the non-compactness of the situation which
requires
careful estimates on $\omega$ given in Lemmas.  In what follows, we
will assume that all solutions of (\ref{41}) are without conjugate
points.

The first step in the proof is the following very well known
construction:  if the solution $u(t)$ has no conjugate points, then
one can easily construct a {\it non-vanishing} solution $\xi$ of the
linearised Jacobi equation
$$
\xi'' + e^{2t} W''_{uu} (u(t),t) \xi = 0
$$
Having such a $\xi(t)$ for every $u(t)$ one defines the function
$\omega (p,u,t)$ by the formula
$$
\omega (p,u,t) = {\dot \xi (t) \over \xi (t)} \quad \mbox{when} \quad
p = \dot u(t), \; u = u(t).
$$
Then $\omega$ satisfies the Ricatti equation along the flow of
(\ref{41}).
\begin{equation}
\dot \omega + \omega^2 + e^{2t} W''_{uu} (u (t), t) = 0
\label{43}
\end{equation}
Here $\cdot$ stands for the derivative along the flow.  It should be
mentioned that by the construction of $\xi$ and $\omega$ the function
$\omega$ is a-priori only measurable (and smooth along the flow).

Denote by
$$
K = \sqrt{sup_{(t,u) \epsilon \Pi} W''_{uu}(u,t)}
$$

We will need the following two lemmas specifying the behaviour of the
function $\omega$ at infinity:

%Lemma 4.4
\begin{lemma}
\label{L44}
The following statements hold true:
%\end{lemma}
% (i), (ii)
\begin{eqnarray}
|\omega (p,u,t)| \leq Ke^T \quad \mbox{for all} \quad (p,u,t)\\
0 \leq \omega (p,u,t) < {1 \over t - T} \quad \mbox{for} \quad t>T
\label{i-ii}
\end{eqnarray}
There exists a constant $\tilde K$ such that for all $(p,u,t)$ with  $t
< T$
% (iii)
\begin{eqnarray}
- \tilde K e^{t/4} < \omega (p,u,t) < Ke^t
\label{iii}
\end{eqnarray}
\end{lemma}
%Lemma 4.5
\begin{lemma}
\label{L45}
The function $\omega$ satisfies the inequalities

I.  For $t \leq T$
\begin{eqnarray}
&&\mbox{if}\quad  u>U, \; p>0 \quad \mbox{then} \quad 0\leq \omega <
  Ke^{t - {u-U \over p}}\label{I1} \\
&&\mbox{if}\quad u<-U, \; p<0 \quad \mbox{then} \quad 0\leq \omega <
  Ke^{t + {-U -u \over p}}\label{I2}\\
&&\mbox{if}\quad u>U -p (T-t), \; p\leq 0 \quad \mbox{or} \quad u< -U-p
  (T-t), \; p \geq 0 \quad \mbox {then} \quad \omega \equiv 0.\label{I3}
\end{eqnarray}
II.  For $t > T$
\begin{eqnarray}
&&\mbox{if}\quad  u>U + p (t-T), \; p>0 \quad \mbox{then} \quad 0\leq
\omega < Ke^{t - {u-U \over p}}\label{II1}\\
&&\mbox{if}\quad  u<U + p (t-T), \; p<0 \quad \mbox{then} \quad 0\leq
\omega < Ke^{t + {-u-U \over p}}\label{II2}\\
&&\mbox{if}\quad  u<-U, \; p \geq 0 \quad \mbox{or} \quad u>U, \; p
\leq 0 \quad \mbox {then} \quad \omega \equiv 0.\label{II3}
\end{eqnarray}
\end{lemma}

We postpone the proof of Lemmas {\ref{L44}} and {\ref{L45}} and first
finish the proof of the theorem.

\noindent{\bf Proof of the theorem}

In order to achieve decay also in the $p$ direction, introduce the
Gibbs density
$$
\alpha (p,u,t) = e^{-H} = e^{- {1/2} {p^2} - e^{2t} W(u,t)}
$$
We have
\begin{equation}
\dot \alpha = -e^{-H} \dot H =  -e^{-H} H_t = - \alpha \left(
  e^{2t} W \right)_t.
\label{3star}
\end{equation}

Now recall that the Hamiltonian flow of (\ref{41}) preserves the
Liouville measure $d \mu = dpdu$.  Multiply the Ricatti equation
(\ref{43}) by $\alpha$ and write it in the form
$$
{\stackrel{\dot \frown}{\alpha \omega}} - \dot \alpha \omega + \alpha
\omega^2 + \alpha e^{2t} W''_{uu} = 0
$$
Substitute equation~(\ref{3star}) to obtain
$$
{\stackrel{\dot \frown}{\alpha \omega}} + \alpha \omega \left( e^{2t}
  W \right)_t + \alpha \omega^2 + \alpha e^{2t} W''_{uu} (u,t) = 0
$$
Owing to the estimates of Lemmas {\ref{L44}} and {\ref{L45}} one can
integrate this equation over the whole $(p,u)$ space.  Use in addition
the invariance of the measure and write
\begin{equation}
{d \over dt} \left(\int \alpha \omega d \mu \right) + \int \alpha
  \omega \left(
  e^{2t} W \right)_t d \mu + \int \alpha \omega^2 d \mu + \int \alpha
  e^{2t} W''_{uu} d \mu = 0
\label{46}
\end{equation}
Using integration by parts the last term can be replaced by $\int
  \alpha \left(e^{2t} W_u \right)^2 d \mu$.

Integrate now the last equation (\ref{46}) for $-A \leq t \leq A$ for
a large constant $A$ and pass to the limit $A \rightarrow + \infty$.
Note that by the uniform estimate of Lemma {\ref{L44}} the term $\int
\alpha \omega d \mu |^A_{-A}$ vanishes in the limit.  So we have
\begin{equation}
\int \alpha \omega \left( W e^{2t} \right)_t d \mu dt + \int \alpha
\omega^2 d \mu dt + \int \alpha \left( e^{2t} W \right)^2 d \mu dt=0
\label{47}
\end{equation}
By the Cauchy--Schwarz inequality we can estimate the first integral of
(\ref{47}) by
$$
\int \alpha \omega \left( W e^{2t} \right)_t d \mu dt \geq - \left[
 \int \alpha \omega^2 d \mu dt \int \alpha \left( \left( We^{2t}
    \right)_t \right)^2 d \mu dt \right]^{1/2}
$$
With the notation $x= \left( \int \alpha \omega^2 d \mu
  dt\right)^{1/2} $ we have the quadratic inequality
$$
x^2 - x \cdot \int \alpha \left[ \left( We^{2t}\right)_t \right]^2 d
\mu dt + \int \left( e^{2t} W_u \right)^2 \alpha d \mu dt \leq 0.
$$
Then its discriminant must be non-negative:
\begin{equation}
4 \int \alpha \left(e^{2t} W_u \right)^2 d \mu dt \leq \int \alpha
\left[ \left( e^{2t} W \right)_t \right]^2 d \mu dt
\label{48}
\end{equation}

The final argument in the proof is the following rescaling
trick.  It is similar to one invented in \cite{[B-P]} for periodic
potentials.
Consider the family of Hamiltonians for every natural number $N$
$$
H_N = {1 \over N^2} H \left( Np, Nu, t \right) = {1 \over 2} p^2 +
{1 \over N^2} e^{2t} W\left( Nu, t \right)
$$
It can be immediately checked that the property of having all
solutions without conjugate points remains valid for all $N$.  Thus
the inequality (\ref{48}) implies the following inequalities for all
$N$:
\begin{eqnarray*}
4 &\int& e^{-{1 \over N^2} W \left( Nu, t \right)e^{2t}} \left(
  e^{2t} {1 \over N} W'_u \left(Nu, t \right) \right)^2 dudt \leq\\
  &\int& e^{-{1 \over N^2} W \left( Nu, t \right)e^{2t}} \left[ {1
  \over N^2} \left( e^{2t} W \left(Nu, t \right) \right)_t\right]^2 dudt
\end{eqnarray*}
(Here the integration with respect to $p$ has been performed on both sides.)

Change the variable in both integrals to $\upsilon=Nu$.  We
obtain the inequality:
\begin{eqnarray*}
4\over N^3 &\int& e^{-{1 \over N^2} W (\upsilon, t) e^{2t}} \left(
  e^{2t} W'_u (\upsilon, t) \right)^2 d\upsilon dt \leq\\
1\over N^5 &\int& e^{-{1 \over N^2} W (\upsilon, t) e^{2t}} \left[ \left(
  e^{2t} W (\upsilon, t) \right)_t \right]^2 d \upsilon dt
\end{eqnarray*}
Now it is clear that if $W$ is not zero identically then the left side
is of order $1 \over N^3$ while the right side is of order $1 \over
N^5$ as $N \rightarrow \infty$.  This proves then that $W \equiv 0$
identically.  The proof of the theorem is completed. $\Box$

%Lemma 4.4
\noindent{\bf Proof of Lemma {\ref{L44}}}

By the definition of $K$ we have
\begin{eqnarray*}
&&\dot \omega \leq K^2 e^{2t} - \omega^2 \quad \mbox{for} \quad t \leq T
\quad \mbox{and}\\
&&\dot \omega = -\omega^2 \quad \mbox{for} \quad t > T
\end{eqnarray*}
The proof is based on the following elementary

\noindent{\bf Fact:}  Any solution of the inequality
$$
\dot \omega \leq B^2 - \omega^2, \quad \omega \left( t_0 \right) =
\omega_0
$$
blows up if $|\omega_0|>B$.  Moreover, the blow up time $t_*$ is
estimated as follows
\begin{eqnarray*}
&&\mbox{for} \quad \omega_0 > B, \quad  t_0 + \Delta < t_* < t_0\\
&&\mbox{for} \quad \omega_0 < -B, \quad  t_0 < t_* < t_0 + \Delta\\
&&\mbox{where} \quad \Delta={1 \over 2B} ln \left( {\omega_0 -B \over
    \omega_0 + B} \right)
\end{eqnarray*}
To prove (i) one takes $B=Ke^T$ and obtains the required estimate.  In the
same manner one obtains (ii) and the right hand side of (iii).

Let us prove the rest of (iii).  Pick two moments of time $\tau_0 <
\tau_1 < 0$.  On the segment $\tau \in \left[ \tau_0, \tau_1
\right]$ we have $\dot \omega \leq B^2 - \omega^2$ with $B=Ke ^{\tau_1}$.  If
$\omega \left( \tau_0 \right) = \omega_0$ is too negative then the
blow up happens before $\tau_1$ and this is impossible.  Thus one
obtains
$$
\Delta = {1 \over 2B} ln \left( {\omega_0 - B \over \omega_0 + B }
\right) > \tau_1 - \tau_0.
$$
This implies
$$
\omega_0 \geq -B {1 + e^{2Bd} \over e^{2Bd} - 1}
$$
Choose $\tau_1 = \tau_0/2$, then this inequality can be written in the
form:
$$
\omega_0 \geq - e^{\tau_0/4} \cdot f \left( \tau_0 \right)
$$
where the function
$$
f \left( \tau_0 \right) = Ke^{\tau_0/4} \cdot {1 + e^{-K \tau_0
    e^{\tau_0/2}} \over 1 + e^{-K \tau_0 e^{\tau_0/2}} }
$$
An easy calculation shows that for $\tau_0 \rightarrow - \infty , \quad
f \rightarrow 0$ and thus $f \left( \tau_0 \right)$ is bounded from
above by some positive constant $\tilde K$.  Since $\tau_0$ was arbitrary,
(iii) follows. $\Box$

%Lemma 4.5
\noindent{\bf Proof of Lemma {\ref{L45}}}

Any point $(p,u,t)$ with $|u|> U$ is situated outside the
support $\Pi$ and so moves in the straight line $u(t) = u+pt$ unless it
hits $\Pi$.  If it does not touch $\Pi$ then $\omega$ vanishes
identically.  This is the case in both (\ref{I3}) and (\ref{II3}).

In the cases of (\ref{I1} -- \ref{I2}) and (\ref{II1} -- \ref{II2})
the line $u(t) = u+pt$ hits $\Pi$ in backward time.  Since for the
free motion $\dot \omega = -\omega^2$, it follows that $\omega (p,u,t)$
has to be non-negative and can be bounded from above by $Ke^{\tilde
  t}$, where $\tilde t$ is the time of entering $\Pi$.

This completes the proof of Lemma {\ref{L45}}. $\Box$

Though by Theorem 4 there are always solutions with conjugate points one cannot
claim that they appear on those radial solutions which are regular
at zero.  The next example shows that all regular solutions may remain
minimal.

%4.6
\noindent{\bf Example {\ref{46}}}

Consider two compactly supported functions $\Phi (u)$ and $\Psi (t)$.
Define the function
$$
W(u,t) = -e^{-2t} \left( \dot \Psi (t) \Phi (u) + {1 \over 2} \Psi (t)
\left( \dot \Phi (u) \right)^2 \right)
$$
Then $W(u,t)$ is compactly supported and one can easily verify that
all the solutions of
\begin{equation}
\dot u = f (u,t) \quad \mbox{with} \quad f (u,t) = \dot \Phi (u) \Psi (t)
\label{last}
\end{equation}
are the solutions of the Newton equation $\mbox{\it \"u} = -e^{2t}
W'_u (u,t)$.

The solutions of (\ref{last}) form a foliation of $\Bbb R(u) \times
\Bbb R(t)$ and then by Weierstrass' theorem of calculus of variations
(see e.g.~\cite{[Hi]} \S 23) are
globally minimal.  It is clear from the construction that all
corresponding radial solutions $u(r)$  of (\ref{11}) are regular at
zero.

\section{Discussion and open problems}
The results of this paper leave open some very natural questions.
We have formulated our results in the context of compactly
supported potentials $V$ on $\Bbb R \times \Bbb R^n$, but most of them
have generalizations to some non-compactly supported cases. Theorem 1
and its proof apply verbatim to $V$ of non-compact support,
in particular to $V$ periodic in $u$, but periodicity in $x$ is
excluded by each of the hypotheses (A) and (B).

Theorem 2 can be generalized to $V$ periodic in u and compactly
supported in $x$, but it is not clear whether it can be extended to other
cases.

It would be very interesting to obtain results for potentials periodic in
both $u$ and $x$ because of the fundamental papers \cite{[M1]}, \cite{[M2]}
by J Moser, which for our equation imply the existence of minimal
laminations for all irrational slopes which for certain slopes are foliations.
In this context our main question looks as follows:  are there
other periodic potentials except those with $V'_u (u,x)=0$ such that for any
slope there is a smooth foliation of $\Bbb T^{n+1}$ by minimal solutions?
We refer the
reader to the survey article by V Bangert \cite{[Ba]} for detailed discussions
of this and related questions.

The
proof of Theorem {\ref{T15}} is based on the reduction to the case of
Hamiltonian systems and so cannot be generalised in a straight-forward
way to non-radial
potentials $V$.  However, it might be that the result remains true.
For example, it is very reasonable to expect that the equation perturbed
by a compactly supported $W$
$$
\Delta u = -V'_u (u,r) + \varepsilon W'_u \left( u, x_1, x_2 \right)
$$
always has non-minimal solutions for $\varepsilon$ sufficiently
small.  In the case of Hamiltonian systems it would follow from the
theorem on continuous dependence of solutions on the initial values.

\hfill

\end{document}